\documentclass[11pt,a4paper]{article}
\usepackage[top=3cm, bottom=3cm, left=3cm, right=3cm]{geometry}
\usepackage[latin1]{inputenc}
\usepackage{amsmath}
\usepackage{amsthm}
\usepackage{amsfonts}
\usepackage{amssymb}
\usepackage{graphics}
\usepackage[hang,flushmargin]{footmisc} 

\usepackage{amsmath,amsthm,amssymb,graphicx, multicol, array}
\usepackage{enumerate}
\usepackage{enumitem}
\usepackage{xcolor}

\usepackage[pagebackref,bookmarks,colorlinks,breaklinks]{hyperref}

\hypersetup{linkcolor=blue,citecolor=blue,filecolor=blue,urlcolor=blue} 

\usepackage{tabto}

\usepackage{tikz}
\usepackage{pgfplots}

\DeclareMathOperator{\li}{li}

\newtheorem{thm}{Theorem}[section]
\newtheorem{lem}{Lemma}[section]
\newtheorem{conj}{Conjecture}[section]

\newtheorem{dfn}{Definition}[section]

\newcommand{\N}{\mathbb{N}}

\newcommand{\R}{\mathbb{R}}
\newcommand{\C}{\mathbb{C}}

\usepackage{fancyhdr}
\usepackage{titletoc}

\titlecontents{section}
[0pt]                                               
{}%
{\contentsmargin{0pt}                               
	\thecontentslabel\enspace%
	\normalsize}
{\contentsmargin{0pt}\normalsize}                        
{\titlerule*[.5pc]{.}\contentspage}                 
[] 
\title{Sum Of Divisors Function Inequality}
\date{}
\author{N. A. Carella}

\begin{document}
\thispagestyle{empty}
\maketitle


\begin{abstract}This short note provides a sharper upper bound of a well known inequality for the sum of divisors function. This is a problem in pure mathematics related to the distribution of prime numbers. 
\let\thefootnote\relax\footnote{ \today \date{} \\
		\textit{AMS MSC2020}: Primary 11A25; Secondary 11A41, 11M06 \\
		\textit{Keywords}: Sum of divisors function, Ramanujan-Robin inequality, Colossally abundant integers, Distribution of prime numbers, Quasi-RH.}
\end{abstract}

\tableofcontents

\section{Introduction} \label{s1}
Let $n \geq 1$ be an integer, and let $\sum_{d \mid n} d$ be the sum of divisors function. The earliest work on the extreme values of the sums of divisors functions appears to be the limit supremum
\begin{equation}
\lim_{n \to \infty} \sup \frac{\sigma(n)}{n \log \log n}=e^{\gamma},
\end{equation}
which was established by Gronwall, see \cite{GT13}, \cite[p.\ 350]{HW79}. Subsequently, conditional on the Riemann hypothesis, Ramanujan proved the sum of divisors function inequality
\begin{equation}
\sigma(n)\leq e^{\gamma}n \log \log n+O((\log n)^{-1/2})
\end{equation}
for all integers $n \geq 3$, see \cite{RS1915}, \cite[Equation 382]{RS1997}. The explicit result, better known as the Ramanujan-Robin inequality, states the same inequality for
 any integer $n\geq 5041$, see \cite{RG84}. Various partial proofs are given in the literature, confer \cite{AJ07}, \cite{BN07}, \cite{BT15}, \cite{CS06}, \cite{GR07},\cite{HA16}, \cite{NS12}, \cite{SP10},  and \cite{WM07}. \\
 
The extreme values of the sum of divisors function occur on a subset of highly composite numbers. The minima occur on the set of primes, and the maxima occur on the set of colossally abundant integers. Extremely abundant integers, colossally abundant integers, etc, are integers related to the primorial integers $n=2^{v_2} \cdot 3^{v_3} \cdots p^{v_p}$, where $p_k$ is the $k$th prime, and $v_k \geq 1$. The precise classifications of the various classes of highly composite integers are given in \cite{RS1915}, \cite{AE44}, \cite{LJ00}, and \cite{BK06}, et alii. Other related results appear in \cite{BN07}, \cite{CS06}, \cite{WM07}. This note proposes an improved result for the upper bound.

\begin{thm} \label{thm1.1A}\hypertarget{thm1.1A}  If $n \geq 3$ is an integer then
	\begin{equation}\label{eq108A}
		\frac{\sigma(n)}{n} \leq e^{\gamma}\log \log n +O\left(\frac{\log\log n}{(\log n)^{1/8} }\right) , 
	\end{equation}	
	where $c>0$ is an absolute constant.
\end{thm}
The proof is given in \hyperlink{S900}{Section} \ref{S900}, it is derived from the new zerofree region of the zeta function, known as the quasi RH \cite{OI2026}. Currently the best unconditional estimates of this arithmetical function are the following.
\begin{thm} \label{thm1.1}\hypertarget{thm1.1}  If $n \geq 3$ is an integer then
\begin{equation}\label{eq108}
\frac{\sigma(n)}{n} \leq e^{\gamma}\log \log n +O\left(e^{-c\sqrt{\log\log n} }\right) , 
\end{equation}	
where $c>0$ is an absolute constant.
\end{thm}
The proof is given in \hyperlink{S900}{Section} \ref{S900}, it is derived from the previous zerofree region of the zeta function, \cite{FK2019}. Another known upper bound is stated below. 
\begin{thm} \label{thm1.2}\hypertarget{thm1.2} Let $n \geq 2521$, then $\sigma(n)/n\leq e^{\gamma}  \log \log n+c/(\log \log n)^2$, 
where $c>0$ is a small constant.
\end{thm}

The above result appears in several papers \cite{RG84}, \cite{NS12}, \cite{AC17} and very recently in \cite{AN2023}, but with different constants. The numerical data for $n \in [5041,  10^{10^{10}}]$ was compiled in \cite{BK06}. Beside the small improvements produced by the small improvements among the zero free regions of the zeta function, see \eqref{eq6363.150k}, the inequality \eqref{eq108A} is the best unconditional result possible. \\

On the other hand, there are several conditional criteria; some of these are listed below. 
\begin{thm} \label{thm1.3}\hypertarget{thm1.3} {\normalfont (\cite{RG84})} Let $n \in \N$ be an integer, and let $\sigma(n)=\sum_{d \mid n} d$ be the sum of divisors function. Then
\begin{enumerate} [font=\normalfont, label=(\roman*)]
\item If the Riemann Hypothesis is true, then, for each $n \geq 5041,$
$$
\sigma(n)< e^{\gamma} n \log \log n.
$$
\item If the Riemann Hypothesis is false, then, there exists constants $0 < \beta < 1/2$ and $c > 0$ such that
$$
\sigma(n)\geq e^{\gamma} n \log \log n+\frac{cn \log \log n}{( \log n)^{\beta}}$$
holds for infinitely many $n$.
\end{enumerate}
\end{thm}

The parameter $1-b < \beta < 1/2$ arises from the possibility of a zero $\rho = b + it \in \C$ of the analytic continuation of
the zeta function 

\begin{equation} \label{4095}
 \zeta(s)=\frac{1}{1-2^{1-s}} \sum_{n \geq 1} \frac{(-1)^{n+1}}{n^s}
\end{equation}  

on the half plane $b = \Re e(\rho) > 1/2$ if the Riemann hypothesis is false. This in turns implies the existence of more or fewer primes than expected in some intervals. For example, under this condition, the number of primes would be 
\begin{equation} \label{4093}
\pi(x)=\li(x)+O(x^{b+\varepsilon})>\li(x)+O(x^{1/2+\varepsilon})
\end{equation}
infinitely often, some of this material is discussed in \cite{LJ00}. The effect of the zeros of the zeta function on the distribution of primes is readily revealed by the explicit formulas, consult the literature. All these inequalities involves the Euler constant $\gamma$. Another well known conditional result offers a constant free inequality, and it is written entirely in terms of the integer $n\geq 1$.

\begin{thm} \label{thm1.4} {\normalfont (\cite{LJ00})}  Let $H_n=\sum_{m \leq n} 1/m$ be the harmonic series. For each $n\geq 1$, the inequality

\begin{equation} \label{4097}
 \sigma(n) < e^{H_n} \log H_n+H_n
\end{equation}

is equivalent to the Riemann Hypothesis.
\end{thm}

The proof of \hyperlink{thm1.1A}{Theorem} \ref{thm1.1A} is given in \hyperlink{S900}{Section} \ref{S900}. The other sections are optional. These sections cover some background materials focusing on the sum of divisors function, and some associated finite sums and products over the prime numbers.

\section{Representations and Identities} \label{S200}\hypertarget{S200}
The sum of divisors function $\sum_{d \mid n}d$ is ubiquitous in number theory. It appears in the analysis many different problems in pure and applied mathematics. Its product representation in (\ref{203}) unearths its intrinsic link to the distribution of the prime numbers. The totient function is defined by $\varphi(n)=\#\{k: \gcd(k,n)=1\}$.

\begin{lem} \label{lem200.1}\hypertarget{lem200.1} Let $n \geq 1$ be an integer, and let the symbol $p^v \mid \mid n$ denotes the maximum prime power divisor. Then
\begin{enumerate} [font=\normalfont, label=(\roman*)]
\item The Euler totient function has the product formula
\begin{equation} \label{201}
\varphi(n)=n\prod_{p \mid n} \left ( 1 -\frac{1}{p} \right ).
\end{equation}
\item The sum of divisors function has the product formula
\begin{equation} \label{203}
\sigma(n)=\prod_{p^v \mid \mid n} \left ( 1 +p+p^{2} +\cdots +p^{v}\right ).\end{equation}
\end{enumerate}
\end{lem}
These representations of the divisor and totient functions are well known and/or easy to establish. 

\begin{lem} \label{lem200.2} Let $n \geq $ be an integer, and let the symbol $p^v \mid \mid n$ denotes the maximum prime power divisor. Then, the sum of divisors function has the product formula
\begin{equation} 
\frac{\sigma(n)}{n}=\frac{n}{\varphi(n)}\prod_{p^v \mid \mid n} \left ( 1 -\frac{1}{p^{v+1}} \right ).
\end{equation}

\begin{proof} Use \hyperlink{lem200.1}{Lemma} \ref{lem200.1}, and the geometric series formula $\sum_{0\leq k \leq x}r=(1-r^{x+1})/(1-r)$ to evaluate the sigma-phi identity
\begin{eqnarray} 
\frac{\sigma(n)}{n}\frac{\varphi(n)}{n}&=&\prod_{p^v \mid \mid n} \left ( 1 +\frac{1}{p}+\frac{1}{p^{2}} +\cdots +\frac{1}{p^{v}}\right ) \prod_{p \mid n} \left ( 1 -\frac{1}{p}\right ) \nonumber \\
&=&\prod_{p^v \mid \mid n} \left ( 1 -\frac{1}{p^{v+1}} \right ) .
\end{eqnarray}
\end{proof}
\end{lem}
To illustrate the negligible effect of negligible multiplication by a prime power, a quantitative expression is computed in the next result.

\begin{lem} \label{lem200.4}\hypertarget{lem200.4} Let $p\geq 2$ be a prime, and let $t\geq 0$ be an integer. If $p^v \mid \mid n$, then
 \begin{equation}
\frac{\sigma(p^{t}n)}{p^{t}n} =\frac{\sigma(n)}{n}  \cdot \frac{1-p^{-(t+v+1)}}{1-p^{-(v+1)}} .
\end{equation}
\end{lem}

\begin{proof}[\textbf{Proof}] By hypothesis $p^v \mid \mid n$. Since $\sigma(n)$ is multiplicative, it is sufficient to observe the effect of extra prime power factor $p^t$ with $t\geq 0$ on the value $\sigma(p^v)/p^v$. To achieve this goal, modify the basic fact $\sigma(p^v)=(p^{v+1}-1)/(p-1)$ to isolate the effect of multiplication by $p^a$. That is,
\begin{eqnarray}
\frac{\sigma(p^{t+v})}{p^{t+v}} &=& \frac{p^{t+v+1}-1}{p^{t+v}(p-1)} \cdot \frac{p^{v+1}-1}{p^{v}(p-1)} \cdot \frac{p^{v}(p-1)}{p^{v+1}-1}   \\
&=& \frac{\sigma(p^{v})}{p^{v}}  \cdot \frac{p^{t+v+1}-1}{p^{t+v}(p-1)} \cdot \frac{p^{v}(p-1)}{p^{v+1}-1}\nonumber \\
&=& \frac{\sigma(p^{v})}{p^{v}}  \cdot \frac{p^{t+v+1}-1}{p^{t}} \cdot \frac{1}{p^{v+1}-1} \nonumber \\
&=& \frac{\sigma(p^{v})}{p^{v}}  \cdot \frac{p^{t+v+1}-1}{p^{t+v+1}-p^t}  \nonumber \\
&=& \frac{\sigma(p^{v})}{p^{v}}  \cdot \frac{1-p^{-(t+v+1)}}{1-p^{-(v+1)}}  \nonumber.
\end{eqnarray}
\end{proof} 

It is immediate that $f(t)=\sigma(p^{t}n)/p^{t}n<2 \log \log n$ is a slowly increasing function of $t\geq 0$, but it is bounded above.


\section{Prime Numbers Theorems} \label{s200}\hyperget{s200}
The omega notation $f (x)= g(x)+\Omega_{\pm}(h(x))$ means that both $f (x) > g(x)+ c_0h(x)$ and $f (x) < g(x)- c_1h(x)$ occur infinitely often as $x \to \infty$, where $c_0 > 0$ and $c_1 > 0$ are constants, see \cite[p.\ 5]{MV2007}, and similar references.\\

The weighted primes counting functions, psi $\psi(x)$ and theta $\theta(x)$, are defined by
\begin{equation}
	\theta(x)=\sum_{p \leq x} \log p \quad \text{ and } \quad \psi(x)=\sum_{p^k \leq x} \log p^k                 
\end{equation}
respectively. The standard prime counting function is denoted by
\begin{equation}
	\pi(x)=\#\{p \leq x\}=\sum_{p \leq x} 1.                 
\end{equation}
This function is usually expressed in term of the logarithm integral $\li(x)=\int_2^x (\log t)^{-1} dt$. The zerofree region Quasi RH 7/8 of the zeta function is proved in \cite{OI2026} and the consequences of for the prime number theorem are discussed in \cite{QRH2026}.

\begin{thm} \label{thm200.1}\hypertarget{thm200.1} Uniformly for $x \geq 2$ the psi and theta functions have the followings asymptotic formulae.
	\begin{enumerate} [font=\normalfont, label=(\roman*)]
		\item  $\displaystyle \theta(x)=x +O\left (x^{7/8}(\log x)^2\right ),$\tabto{9cm} Unconditionally,	
		\item $\displaystyle\theta(x)=x +\Omega_{\pm} \left ( x^{1/2}\log \log \log x \right ),$\tabto{9cm}Unconditional oscillation,
		\item
		$\displaystyle\theta(x)=x +O\left (x^{1/2} \log^2 x \right ),$\tabto{9cm} Conditional on the RH.
	\end{enumerate}
\end{thm}

\begin{proof} (i) The new zerofree region $\{s\in\C:\Re e(s)>7/8\}$ of the zeta function $\zeta(s)$, see \cite[Theorem 1.1]{OI2026} and an effective version of the explicit formula immediately yield this result.\\
	
(ii) The oscillations form of the theta function is proved in \cite[p.\ 479]{MV2007}.

(iii)The proof of this case is similar to the first case (i).
\end{proof}

The same asymptotics hold for the function $\psi(x)$. Explicit estimates for both of these functions are given in \cite{CP1985}, \cite{SL1976}, \cite[Theorem 5.2]{DP2010}, and related literature. 

\begin{conj} \label{conj200.1} Assuming the RH and the LI conjecture, the suprema are
	\begin{equation}
		\lim \inf_{x \to \infty} \frac{\psi(x)-x}{\sqrt{x} (\log \log x)^2}=\frac{-1}{\pi} \qquad \text{and} \qquad \lim \sup_{x \to \infty} \frac{\psi(x)-x}{\sqrt{x} (\log \log x)^2}=\frac{1}{\pi}.
	\end{equation}
\end{conj}

More details on the Linear Independence conjecture appear in \cite{IA1942}, \cite[Theorem 6.4]{EE1985}, and recent literature. The LI conjecture asserts that the imaginary parts of the nontrivial zeros $\rho_n=1/2+i \gamma_n$ of the zeta function $\zeta(s)$ are linearly independent over the set $\{-1,0, 1\}$. In short, the equations  
\begin{equation}
	\sum_{1 \leq n \leq M} r_n \gamma_n=0,
\end{equation}
where $r_n \in \{-1,0, 1\}$, have no nontrivial solutions. 

\begin{thm} \label{thm200.2}\hypertarget{thm200.2} Let $x \geq 1$ be a large number. Then
	\begin{enumerate} [font=\normalfont, label=(\roman*)]
		\item  $\displaystyle \pi(x)=\frac{x}{\log x} +O\left (x^{7/8}\log x\right ),$\tabto{9cm} Unconditionally,
		\item $\displaystyle\pi(x)=\li(x) +\Omega_{\pm} \left (\frac{ x^{1/2}\log \log \log x}{\log x} \right ),$\tabto{9cm}Unconditional oscillation,
		
		\item
		$\displaystyle\pi(x)=\li(x) +O\left (x^{1/2} \log x \right ),$\tabto{9cm} Conditional on the RH.
	\end{enumerate}
\end{thm}

\begin{proof}[\textbf{Proof}] (i) This follows from an application of \hyperlink{thm200.1}{Theorem} \ref{thm200.1} and summation by part. The earlier unconditional prime counting function $\pi(x)=\li(x)+O\left ( x e^{-c_0 \log x^{3/5}(\log \log x)^{-2/5}}\right )$ appears in \cite{FK2019}. \\
	
	(ii) The unconditional oscillations part arises 
	from the Littlewood form $\pi(x)=\li(x)+\Omega_{\pm} \left (x^{1/2}\log \log \log  x /\log x \right )$ of the prime number theorem, 
	consult \cite[p.\ 51]{IA2003}, \cite[p.\ 479]{MV2007}, et cetera. \\
	
	(iii)  The conditional part arises from the Riemann form $\pi(x)=\li(x)+O\left (x^{1/2}\log^2 x \right )$ of the prime number theorem. 
\end{proof}

New explicit estimates for the number of primes in arithmetic progressions are computed in \cite{BR2018}.

\section{Harmonic Sums over the Primes} \label{s2} \hyperget{s2}
The most basic finite sum over the prime numbers is the prime harmonic sum $\sum_{p \leq x}1/p$. The refined estimate of this finite sum, stated below, is a synthesis of various results due to various authors. 

\begin{lem} \label{lem2.1}\hyperget{lem2.1}
	Let \(x\geq 2\) be a large number, then
	\begin{enumerate} [font=\normalfont, label=(\roman*)]
		\item  $\displaystyle\sum_{p \leq x} \frac{1}{p}
		=\log \log x+ B_1 +O\left (\frac{\log x}{x^{1/8}} \right ), $ \tabto{8cm}  Unconditionally,
		\item $\displaystyle \sum_{p \leq x} \frac{1}{p}
		=\log \log x+ B_1 +\Omega_{\pm} \left (\frac{\log \log \log x}{x^{1/2}\log x} \right ),  $\tabto{8cm} Unconditional oscillation
		\item $\displaystyle \sum_{p \leq x} \frac{1}{p}
		=\log \log x+ B_1 +O\left (\frac{\log x}{ x^{1/2}} \right ),$\tabto{8cm}Conditional on the RH,
	\end{enumerate}
	where $B_1 = 0.2614972128 \ldots, $ is Mertens constant, and $c_0>0$ is an absolute constant. 
\end{lem}

\begin{proof}[\textbf{Proof}] (i) Replace prime counting measure $\pi(x)=\li(x)+O(x^{7/8}\log x)$, in \hyperlink{thm200.2}{Theorem} \ref{thm200.2}, into the Stieltjes integral representation and evaluate it.
	\begin{align}
		\sum_{p \leq x}  \frac{1}{p}&=\int_2^{x}\frac{1}{t}d\pi(t) \\[.3cm]
		&=\frac{\li(x)+O(x^{7/8}\log x)}{x}+\int_2^{x}\frac{\li(t)+O(t^{7/8}\log t)}{t^2}dt \\[.3cm]
		&=\log\log x+B_1+O\left(\frac{\log x}{x^{1/8}} \right), 
	\end{align}
	where $B_1>0$ is a constant.\\
	
	(ii) The unconditional oscillations part arises 
	from the Littlewood form $\pi(x)=\li(x)+\Omega_{\pm} \left (x^{1/2}\log \log \log  x /\log x \right )$ of the prime number theorem, consult \cite[p.\ 51]{IA03}, \cite[p.\ 479]{MV07}, et cetera. \\
	
	(iii)  The conditional part arises from the Riemann form $\pi(x)=\li(x)+O\left (x^{1/2}\log^2 x \right )$ of the prime number theorem. 
\end{proof}

The asymptotic order $\sum_{p \leq x}1/p \sim \log \log x$ is due to Euler, confer \cite[Chapter 15?]{EL1888}. The earliest version including error term $\sum_{p \leq x}1/p=\log \log x+ B_1 +O(1/ \log x)$ is due to Mertens, see \cite{VM2005}. The qualitative form of the oscillations of the differences 
\begin{equation} \label{eq200.8}
	\sum_{p^k \leq x} \frac{1}{p^k}-\left ( \log \log x + \gamma \right ) \quad \text{ and} \quad \sum_{p^k \leq x} \frac{\log p}{p^k}-\left ( \log x + \gamma \right )
\end{equation} 
seems to be due to Phragmen, confer \cite[p.\ 182]{NW2000}. \\

The Euler constant and Mertens constant occur very frequently in analysis. The former is defined by 
\begin{equation}
	\gamma=\lim_{x \to \infty} \left ( \sum_{n \leq x}\frac{1}{n}-\log x \right ) =0.577215665 \ldots,
\end{equation}
and the later 
is defined by 
\begin{equation}
	B_1=\lim_{x \to \infty} \left ( \sum_{p \leq x}\frac{1}{p-1}-\log \log x \right ) =0.2614972128 \ldots.
\end{equation}
Other definitions of these constants are available in the literature, confer \cite{LJ13}.

\begin{lem} \label{lem2.2}\hypertarget{lem2.2} The constants $\gamma$ and $B_1$ satisfy the linear relation
	\begin{equation}
		B_1=\gamma-\sum_{p \geq 2} \sum_{n \geq 2} \frac{1}{np^{n}}.                 
	\end{equation}
	
\end{lem}

\begin{proof}[\textbf{Proof}] This relation stems from the power series expansion 
	\begin{equation}
		B_1-\gamma=\sum_{p \geq 2} \left (\log \left (1-\frac{1}{p} \right ) + \frac{1}{p} \right )
	\end{equation}
	via the power series for $\log(1+z)$ with $|z|<1$. This leads to this identity, see \cite[p.\ 466]{HW79}, \cite[p.\ 182]{MV2007}.\end{proof}

\begin{lem} \label{lem2.3} Let \(x\geq 2\) be a large number, then
	\begin{equation}
		\sum_{p \leq x} \sum_{n \geq 2} \frac{1}{np^{n}} =\gamma-B_1-\frac{1}{x \log x}+ O \left (\frac{1}{x\log^2 x}\right ). 
	\end{equation}
\end{lem}

\begin{proof}[\textbf{Proof}] Rearrange the power series expansion as
	\begin{equation}
		\sum_{p \leq x} \sum_{n \geq 2} \frac{1}{np^{n}} =\gamma-B_1- \sum_{p > x} \sum_{n \geq 2} \frac{1}{np^{n}} =\gamma-B_1-\frac{1}{x \log x}+ O \left (\frac{1}{x\log^2 x}\right ). 
	\end{equation}
	The estimate for the last two terms on the right follows from \hyperlink{lem2.4}{Lemma} \ref{lem2.4} computed below.
\end{proof} 

\begin{lem} \label{lem2.4}  Let $x \geq 1$ be a large number, then
	\begin{equation}
		\sum_{p >x}\sum_{n \geq 2} 
		\frac{1}{np^{n}} = \frac{1}{x \log x}+ O \left (\frac{1}{x\log^2 x}\right ).
	\end{equation}
\end{lem}
\begin{proof}[\textbf{Proof}] Split the infinite sum into two subsums:
	\begin{eqnarray}
		\sum_{p > x} \sum_{n \geq 2} \frac{1}{np^n}&=&\sum_{p > x}  \frac{1}{2p^2}+\sum_{p > x} \sum_{n \geq 3} \frac{1}{np^n} \nonumber \\
		&=&\sum_{p \geq x}  \frac{1}{2p^2}+O \left (\frac{1}{x^2 \log x}\right) \nonumber.
	\end{eqnarray}
	Employ the prime counting measure $\pi(t)=\#\{p \leq t\}$ to evaluate the first subsum using the integral 
	\begin{eqnarray}
		\sum_{p \geq x}  \frac{1}{p^2}&=&\int_x^{\infty}\frac{1}{t^2}d\pi(t) \nonumber \\
		&=&-\frac{\pi(x)}{x^2} +2\int_x^{\infty}\frac{\pi(t)}{t^3}dt \nonumber \\
		&=&\frac{1}{x \log x}+O \left (\frac{1}{x \log^2 x}\right) \nonumber .
	\end{eqnarray} 
\end{proof}

\section{Products over the Primes } \label{s400}
The prime harmonic products are sine qua non in number theory. It has a natural link to the totient function $\varphi(n)$, the sieve of Eratosthenes, and other related concepts. Accurate estimates of this product and related products are essential in a variety of calculations in number theory.

\begin{lem} \label{lem400.1}\hypertarget{lem400.1} {\normalfont (\cite[Theorems 7, 8]{RS1962})}  Let $x \geq 2$ be a large number. Then
\begin{enumerate} [font=\normalfont, label=(\roman*)]
\item 
$\displaystyle\left | \frac{1}{e^{\gamma} \log x}\prod_{p \leq x} \left ( 1 -\frac{1}{p} \right )^{-1} - 1  \right | <  \frac{1}{2( \log x)^2 }.$\\

\item $\displaystyle\left | (e^{\gamma} \log x )\prod_{p \leq x} \left ( 1 -\frac{1}{p} \right ) - 1 \right | <\frac{1}{2( \log x)^2 }.$
\end{enumerate}
\end{lem}

The new zerofree region of the zeta function, known as the quasi RH, \cite{OI2026}, automatically leads to sharper results for the associated sums and products of the prime numbers.
\begin{lem} \label{lem400.1A}\hypertarget{lem400.1A}   Let $x \geq 2$ be a large number. Then
	\begin{enumerate} [font=\normalfont, label=(\roman*)]
		\item 
		$\displaystyle\left | \frac{1}{e^{\gamma} \log x}\prod_{p \leq x} \left ( 1 -\frac{1}{p} \right )^{-1} - 1  \right | \ll  \frac{ \log x }{x^{1/8}}.$\\
		
		\item $\displaystyle\left | (e^{\gamma} \log x )\prod_{p \leq x} \left ( 1 -\frac{1}{p} \right ) - 1 \right | \ll  \frac{ \log x }{x^{1/8}}.$
	\end{enumerate}
\end{lem}
\begin{proof}  (i) The logarithm of the product has the from
\begin{equation}\label{eq400.1A02}
	\log \prod_{p \leq x} \left ( 1 -\frac{1}{p} \right )^{-1}=\sum_{p \leq x} \log \left ( 1 -\frac{1}{p} \right )^{-1}=\sum_{p \leq x}\frac{1}{p}+\sum_{p \leq x}\sum_{k\geq2}\frac{(-1)^{k+1}}{kp^k}.
\end{equation}	
An asymptotic formula for the first finite is compute below and the double sum above converges to a constant. \\	
	
Replace prime counting measure $\pi(x)=\li(x)+O(x^{7/8}\log x)$ into the Stieltjes integral representation and evaluate it.
	\begin{align}
		\sum_{p \leq x}  \frac{1}{p}&=\int_2^{x}\frac{1}{t}d\pi(t) \\[.3cm]
		&=\frac{\li(x)+O(x^{7/8}\log x)}{x}+\int_2^{x}\frac{\li(t)+O(t^{7/8}\log t}{t^2}dt \\[.3cm]
		&=\log\log x+B_1+O\left(\frac{\log x}{x^{1/8}} \right), 
	\end{align}
	where $B_1>0$ is a constant. Reversing the logarithm in \eqref{eq400.1A02} yields the \begin{equation}\label{eq400.1A06}
	 \prod_{p \leq x} \left ( 1 -\frac{1}{p} \right )^{-1}=e^{\gamma} \log x+O\left ( \frac{\log x}{x^{1/8}}\right) .
	\end{equation}	
The is a routine procedure to resolve the constant as $\gamma$, see \hyperlink{lem2.2}{Lemma} \ref{lem2.2}.
\end{proof}
These are improved versions of the original works by Mertens, see \cite{RS1962}, and \cite{VM2005} for more details. The Euler constant is defined by $\gamma= \lim_{x \to \infty} \left ( \sum_{n \leq x} 1/n- \log x \right )=.577215665 \ldots$, see \cite[p.\ 28]{FS03}, \cite[Theorem 2.2.1]{LJ13} for several definitions of this number and similar references.

\begin{lem} \label{lem400.2}\hypertarget{lem400.2} {\normalfont (\cite[Corollary 3]{SL1976})}  Let $x \geq 14$ be a large number. If the Riemann hypothesis holds, then
\begin{enumerate} [font=\normalfont, label=(\roman*)]
\item $\displaystyle 
\left | \frac{1}{e^{\gamma} \log x}\prod_{p \leq x} \left ( 1 -\frac{1}{p} \right )^{-1} - 1  \right | <\frac{ \log x+5}{8 \pi \sqrt{x} }.$\\

\item $\displaystyle \left | (e^{\gamma} \log x )\prod_{p \leq x} \left ( 1 -\frac{1}{p} \right ) - 1 \right | <\frac{3 \log x+5 }{8 \pi \sqrt{x} }.$
\end{enumerate}
\end{lem}

\begin{thm} \label{thm400.3} {\normalfont (Mertens)} Let $n \geq 1$ be an integer. Then
\begin{enumerate} [font=\normalfont, label=(\roman*)]
\item $\displaystyle\lim_{x \to \infty} \sup \frac{1}{ \log x} \prod_{p \leq x} \left ( 1 -\frac{1}{p} \right )^{-1}=e^{\gamma},$\tabto{7cm}The limit supremum over the prime.\\

\item $\displaystyle\lim_{x \to \infty} \sup \frac{1}{ \log x} \prod_{p \leq x} \left ( 1 +\frac{1}{p} \right )^{-1}=\frac{6e^{\gamma}}{\pi^2},$\tabto{7cm}The limit supremum over the prime.
\end{enumerate}
\end{thm}
Some references and other information on these limits are available in  see \cite[p.\ 31]{FS03}. 
The results for products over arithmetic progression are proved in \cite{LZ2007}, et alii.\\

\begin{lem} \label{lem2003.701}
	Let \(x\geq 2\) be a large number, then
	\begin{enumerate} [font=\normalfont, label=(\roman*)]
		\item  $\displaystyle \prod_{p \leq x}\left( 1+ \frac{1}{p} \right) 
		=\frac{6e^{\gamma}}{\pi^2} \log x+ O\left (\frac{\log x}{x^{1/8}}\right ) , $ \tabto{9cm}  unconditionally.
		\item 
		$\displaystyle
		\prod_{p \leq x}\left( 1+ \frac{1}{p} \right) 
		=\frac{6e^{\gamma}}{\pi^2} \log x+\Omega \left (\frac{\log \log \log x}{x^{1/2}\log x} \right ),$\tabto{9cm} unconditional oscillation.
		\item $\displaystyle
		\prod_{p \leq x}\left( 1+ \frac{1}{p} \right) 
		=\frac{6e^{\gamma}}{\pi^2} \log x+ O\left (\frac{\log x}{ x^{1/2}} \right ) ,  
		$\tabto{9cm} conditional on the RH.
	\end{enumerate}
	where $\gamma$ is Euler constant, and $c_0>0$ is an absolute constant. 
\end{lem}
\begin{proof} (i) For a large real number $x \in \R$, rewrite the product as
	\begin{eqnarray}\label{eq2003.701}
		\prod_{p \leq x}\left( 1+ \frac{1}{p} \right)&=& \prod_{p \leq x}\left( 1+ \frac{1}{p} \right)\left( 1- \frac{1}{p} \right) \left( 1- \frac{1}{p} \right)^{-1} \\
		&=& \prod_{p \leq x}\left( 1- \frac{1}{p^2} \right)\prod_{p \leq x}\left( 1- \frac{1}{p} \right)^{-1}\nonumber.
	\end{eqnarray}
	Replacing the completed product
	\begin{equation}\label{eq2003.705}
		\prod_{p \leq x}\left( 1- \frac{1}{p^2} \right)=\prod_{p \geq 2}\left( 1- \frac{1}{p^2} \right)+O\left ( \frac{1}{x^c}\right )=\frac{6}{\pi^2}+O\left ( \frac{1}{x^c}\right ) ,
	\end{equation}
	where $c\geq 1$ is a constant, in the first product on the right side of \eqref{eq2003.701}, and applying \hyperlink{lem400.1A}{Lemma} \ref{lem400.1A} yield
	\begin{eqnarray}\label{eq2003.707}
		\prod_{p \leq x}\left( 1- \frac{1}{p^2} \right)\prod_{p \leq x}\left( 1- \frac{1}{p} \right)^{-1}&=& \left ( \frac{6}{\pi^2}+O\left ( \frac{1}{x^c}\right )\right ) \prod_{p \leq x}\left( 1- \frac{1}{p} \right)^{-1}\\
		&=& \left ( \frac{6}{\pi^2}+O\left ( \frac{1}{x^c}\right )\right ) \left( e^{\gamma} \log x+ O\left (\frac{\log x}{x^{1/8}}\right )   \right) \nonumber \\
		&=& \frac{6e^{\gamma}}{\pi^2} \log x+ O\left (\frac{\log x}{x^{1/8}}\right )   \nonumber.
	\end{eqnarray}
	The verification of statements (ii) and (iii) are similar, mutatis mutandus. 
	
\end{proof}
The nonquantitative unconditional oscillations of the error of the product of primes is implied by the work of Phragmen, refer to equation (\ref{eq3.10}), and \cite[p.\ 182]{NW2000}. Since then, various authors have developed quantitative versions, see \cite{RS1962}, \cite{DP2009}, \cite{LY2014}, \cite{LJ2015}, et alii. The specific quantitative form
\begin{equation} \label{eq3.10}
	\prod_{p \leq x}\left( 1- \frac{1}{p} \right) ^{-1}
	=e^{\gamma} \log x+ \Omega_{\pm}\left (\frac{f(x)}{ x^{1/2}} \right ), 
\end{equation}
where $f(x)$ is a slowly increasing function, was proved in \cite{DP2009}. \\

\begin{thm} \label{thm2003.2} { \normalfont (Martens, 1874)} The following asymptotic formulas hold:
	\begin{enumerate} [font=\normalfont, label=(\roman*)]
		\item  $\displaystyle \lim_{x \to \infty} \frac{1}{\log x} \prod_{p \leq x} \left (1-\frac{1}{p}\right )^{-1}=e^{\gamma}, $ \tabto{7cm}  the product over the integers.
		\item 
		$\displaystyle
		\lim_{x \to \infty} \frac{1}{\log x} \prod_{p \leq x} \left (1+\frac{1}{p}\right )^{-1}=\frac{6e^{\gamma}}{\pi^2}$\tabto{7cm} the product over the squarefree integers.
	\end{enumerate}
\end{thm}

\section{Auxiliary Results}\label{S6363}\hypertarget{S6363}
For a complex number $s\in \C$, and a real number $x\geq1$ the  \textit{rho series} is defined by
\begin{equation}
	S_s(x)=	 -s\sum_{\rho}\frac{x^{\rho-s}}{\rho(\rho-s)},
\end{equation}
where the index $\rho\in \C$ ranges over the nontrivial zeros of the zeta function, see \cite[p.\ 134]{RS1997}.  
\begin{lem}\label{lem6363.200}\hypertarget{lem6363.200} Let $x>1$ be a large number and let $\rho$ denotes the nontrivial zero of the zeta function, then
	\begin{enumerate} [font=\normalfont, label=(\roman*)]
		\item$\displaystyle  \sum_{\rho}\frac{x^{\rho-1}}{\rho(1-\rho)}\ll e^{-c\sqrt{\log x}} ,$ \tabto{8cm} Unconditionally,
		\item$\displaystyle \sum_{\rho}\frac{x^{\rho-1}}{|\rho|^2}\ll\frac{ 1}{\sqrt{ x}},$ \tabto{8cm} Conditioned on the RH,
	\end{enumerate}
	where $c>0$ is a constant.
\end{lem}
\begin{proof}[\textbf{Proof}] (i) Let $\rho=\sigma+it=1-c_0/\log t+it$ be a nontrivial zero of the zeta function, where $t\in [14,T]$, see \cite[p.\ 87]{DH1980}, \cite{IA2003},  \cite{FK2019}, et cetera. For any small number $\varepsilon>0$, let $T = e^{c_1(\log x)^{\varepsilon}}$. Then
	\begin{equation}\label{eq6363.200c}
		x^{-c_0/\log T}\ll e^{-c(\log x)^{1-\varepsilon}}.
	\end{equation}		
	The choice of $\varepsilon=1/2$ leads to the standard error term
	\begin{equation}\label{eq6363.200d}
		x^{-c_0/\log T}\ll e^{-c\sqrt{\log x}},
	\end{equation}	
	where $c_0,c_1,  c>0$ are constants, see the calculation of error term of the prime number theorem in \cite[Theorem 6.9, page 181]{MV2007}. Unconditionally, the series over the nontrivial zeros of the zeta function has the value
	\begin{equation}\label{eq6363.200e}
		B=\sum_{\rho}\frac{1}{\rho(1-\rho)}\ll 1
	\end{equation}
	and conditioned on the RH, the series over the nontrivial zeros of the zeta function has the value
	\begin{equation}\label{eq6363.200f} 
		B=\sum_{\rho}\frac{1}{|\rho|^2}=-\frac{\gamma}{2}-1+\frac{1}{2}\log 4\pi
	\end{equation}
	see \cite[p.\ 163]{DH1980}. Multiplying \eqref{eq6363.200d} and \eqref{eq6363.200e} complete the verification. The verification of statement (ii) is similar to the proof for statement (i) but it uses $\rho=\sigma+it=1/2+it$, where $t\in \R$ and \eqref{eq6363.200f} .
\end{proof}
The sharper estimate
\begin{equation}\label{eq6363.150k} 
	\rho=\sigma+it=1-\frac{c_0}{(\log\log t)^{2/3}(\log\log\log t)^{1/3}}+it
\end{equation}
of the nontrivial zero of the zeta function provides some improvement on these estimates, but these sharper estimates are not required in this note, see \cite[Theorem 1]{FK2019} for the precise details.\\

The earliest conditional result seems to be the well known upper bound of the sum of divisor function shown below.

\begin{thm}\label{thm6363.100}\hypertarget{thm6363.100}{\normalfont (\cite[Equation 382, p.\ 143]{RS1997})} If $n\geq3$ is an integer then
	$$\frac{\sigma(n)}{n}\leq e^{\gamma}\left( \log \log n+\frac{a_0}{\sqrt{\log n}}+S_1(\log n)+\frac{a_1}{\sqrt{\log n}\log \log n}\right),  $$
	where
	$$S_1(x)=\sum_{\rho}\frac{x^{\rho-1}}{|\rho|^2},$$
	and $a_0<0$ and $a_1$ are constants.
\end{thm}

\begin{thm}\label{thm6363.100RH}\hypertarget{thm6363.100RH}{\normalfont (\cite[Equation 382, p.\ 143]{RS1997})} Assume the RH. If $n\geq3$ is an integer then
	$$\frac{\sigma(n)}{n}\leq e^{\gamma}\log \log n+\frac{c}{\sqrt{ \log n}},  $$
	where $c>0$ is a constant.
\end{thm}

\begin{proof} This follows from Lemma \ref{lem6363.200} and Theorem \ref{thm6363.100}.
\end{proof}
\section{Highly Composite Numbers} \label{s500}
Let $p_k$ be the $k$th prime in increasing order, and let $v_p=\max \{ m: p^m \mid n\}$ is the $p$-adic valuation. Extremely abundant integers, and colossally abundant integers are related to the primorial integers $n=2^{v_2} \cdot 3^{v_3} \cdots p^{v_p}$, but the exponents have certain multiplicative structure $1 \leq v_p\leq \cdots \leq v_3 \leq v_2$. 

\begin{dfn} \label{dfn500.4} {\normalfont Let $d(n)=\sum_{d \mid n}1$. An integer $n \in \N$ is called highly composite if and only if $ d(m) 1 < d(n)$ for all integers $m <n$.}
\end{dfn}

\begin{dfn} \label{dfn500.5} {\normalfont Let $\sigma(n)=\sum_{d \mid n}d$. An integer $n \in \N$ is called colossally abundant if and only if
\begin{equation} \frac{\sigma(m)}{m^{1+\varepsilon} }< \frac{\sigma(n)}{n^{1+\varepsilon}}
\end{equation}
for all integers $m <n$, and some small number $\varepsilon>0$}
\end{dfn}
\begin{dfn} \label{dfn500.6} {\normalfont An integer $n \geq 10080$ is called extremely abundant if and only if
\begin{equation} \frac{\sigma(m)}{m \log \log m} < \frac{\sigma(n)}{n \log \log n}
\end{equation}
for all integers $m <n$}
\end{dfn}
These numbers are studied in \cite{RS1915}, \cite{AE44}, \cite{LJ00}, \cite{BK06}, \cite{NS12}, et alii. \\

Fixed a highly composite integer $n \geq 1$, and let $n_r=n^r$ with a parameter $r \geq 1$. The sequence of functions
\begin{equation} \label{eq500.7}
\frac{\sigma(n_r)}{n_r} < \frac{\sigma(n_{r+1})}{n_{r+1}}<\frac{\sigma(n_{r+1})}{n_{r+1}} < \cdots
\end{equation}
is strictly monotonically increasing, but bounded above by $e^{\gamma} \log \log n$, see equation (\ref{eq700.6}). 

\section{Prime Factors of Highly Composite Integers}
\begin{lem} \label{lem500.1} {\normalfont (\cite[Theorem 2]{AE44})} Let $n \geq $ be a large highly composite integer, then
\begin{enumerate} [font=\normalfont, label=(\roman*)]
\item Unconditionally, the largest prime divisor $p \mid n$ has the asymptotic 
$$ p= (\log n)\left (1+O\left (\frac{1}{(\log \log n)^2} \right ) \right ).$$
\item  Modulo the Riemann hypothesis, the largest prime divisor $p \mid n$ has the asymptotic 
$$ p= (\log n)\left (1+O\left (\frac{\log \log n}{(\log n)^{1/2}} \right ) \right ).$$
\end{enumerate}
\end{lem}

\begin{lem} \label{lem500.100} Let $n \geq 1$ be a large highly composite integer, then
\begin{equation}\label{eq500.100}
\frac{\sigma(n)}{n}<e^{\gamma}\log (x) \left(1+\frac{.5}{(\log x)^2}\right ),
\end{equation}
where $x\ll \log n$.
\end{lem}
\begin{proof}[\textbf{Proof}] Let $v=v_p\geq0$ be the valuation of $n$ at $p$. Suppose that each prime divisor $p\mid n$ satisfies $p\leq x$. 
\begin{eqnarray}\label{eq500.110}
\frac{\sigma(n)}{n}&=&\prod_{p^v\mid \mid n}\left(\frac{1-p^{-v-1}}{1-p^{-1}}\right)\\
&=&\prod_{p^v\mid \mid n}\left(1+\frac{1}{p}+\frac{1}{p^2}+\cdots+\frac{1}{p^v}\right)\nonumber\\
&<&\prod_{p\leq x}\left(1-\frac{1}{p}\right)^{-1}\nonumber.
\end{eqnarray}
By \hyperlink{lem400.1}{Lemma} \ref{lem400.1}, it follows that
\begin{equation}\label{eq500.102}
\frac{\sigma(n)}{n}<e^{\gamma}\log (x) \left(1+\frac{.5}{(\log x)^2}\right ),
\end{equation}
where $x\ll \log n$.
\end{proof}
\begin{lem} \label{lem500.500} Let $n \geq 1$ be a large highly composite integer, then
\begin{equation}\label{eq500.500}
\log \log n >\log(x) \left(1+O \left ( \frac{1}{(\log x)^{3}}\right )\right),
\end{equation}
where $x\ll \log n$.
\end{lem}
\begin{proof}[\textbf{Proof}] Assume that $n=\prod_{p^{v_p}\mid \mid n}2^{v_2}\cdot 3^{v_3}\cdots p^{v_p}$.Clearly, for any highly composite integer $v_2\geq 2$. Hence,

\begin{eqnarray}\label{eq500.510}
\log n&=&\log \prod_{p^{v_p}\mid \mid n}2^{v_2}\cdot 3^{v_3}\cdots p^{v_p}\\
&=&\sum_{p^{v_p}\mid \mid n}\log p^{v_p}\nonumber\\
&>&(v_2-1)\log 2+\sum_{p\mid n}\log p\nonumber\\
&\geq&\log 2+\sum_{p\leq x}\log p\nonumber,
\end{eqnarray}
where $x\ll \log n$. By the prime number theorem 
\begin{equation}\label{eq500.520}
\sum_{p\leq x}\log p=x+O\left (\frac{x}{(\log x)^{2}}\right )\geq 
x-\left (\frac{c_0x}{(\log x)^{2}}\right ),
\end{equation}
where $c_0>0$ is a constant. Thus,

\begin{eqnarray}\label{eq500.530}
\log \log n&>&\log \left(\log 2+x+O(x(\log x)^{-2})\right)\\
&\geq&\log \left(\log 2+x-(c_0x(\log x)^{-2})\right)\nonumber\\
&=&\log \left(x \left ( 1+\frac{\log 2-(c_0x(\log x)^{-2})}{x}\right )\right)\nonumber\\
&=&\log(x) \left(1+\frac{\log \left ( 1+\frac{\log 2-(c_0x(\log x)^{-2})}{x}\right )}{\log x}\right)\nonumber\\
&=&\log(x) \left(1- \left ( \frac{c_1}{(\log x)^{3}}\right )\right)\nonumber,
\end{eqnarray}
where $c_1>0$ is a constant.
\end{proof}

\section{The Main Results} \label{S900}\hypertarget{S900}
Let $N_k=2\cdot3\cdots p_k$ be a primorial number, it is the product of the first $k$ primes $p_k$ in increasing order. The theta function has the asymptotic formula $\theta(x)=\sum_{p\leq x}\log p=x+o(x)$. At $x=N_k$, it has the form
\begin{equation}\label{eq7.000}
	\theta(p_k)=\sum_{p\leq p_k}\log p=\log N_k.
\end{equation} 
\begin{proof}[\textbf{Proof of {\normalfont \hyperlink{thm1.1A}{Theorem} \ref{thm1.1A}}}] Take the sigma-phi representation of the sum of divisors function, see  \hyperlink{lem200.2}{Lemma} \ref{lem200.2}, and utilize \hyperlink{lem400.1}{Lemma} \ref{lem400.1} to evaluate the product representation.
	\begin{align}\label{eq900.200}
		\frac{\sigma(n)}{n}&<			\frac{n}{\varphi(n)}\\[.3cm]&=\prod_{p\mid n}\left(1-\frac{1}{p} \right)^{-1}\nonumber\\[.3cm]
		&\leq\prod_{p\mid N_k}\left(1-\frac{1}{p} \right)^{-1}\nonumber\\[.3cm]
		&\leq\prod_{p\ll \log N_k}\left(1-\frac{1}{p} \right)^{-1}\nonumber\\[.3cm]
		&=  e^{\gamma} \log \log N_k+O\left( \frac{\log \log N_k}{( \log N_k)^{1/8}}\right) \nonumber.
	\end{align}
	This completes the verification.
\end{proof}

\begin{proof}[\textbf{Proof of {\normalfont \hyperlink{thm1.1}{Theorem} \ref{thm1.1}}}] Let $n \geq 5041$ be an integer and let $x\asymp\log n$. Substituting the unconditional estimate for the rho series 
	\begin{equation}\label{eq6363.150f}
		S_1(\log n)\ll e^{-c\sqrt{\log \log n}},
	\end{equation}	
	confer \hyperlink{lem6363.200}{Lemma} \ref{lem6363.200}, into the expression given in \hyperlink{thm6363.100}{Theorem} \ref{thm6363.100} yields
	\begin{eqnarray}\label{eq6363.150e}
		\frac{\sigma(n)}{n}	&\leq&e^{\gamma}\left( \log \log n+\frac{a_0}{\sqrt{ \log n}}+S_1(\log n)+\frac{a_1}{\sqrt{\log n}\log \log n}\right)\nonumber\\		
		&=&e^{\gamma} \log \log n+O\left(e^{-c\sqrt{\log \log n}}\right)
	\end{eqnarray}	
	where $c>0$ is a constant.
\end{proof}
As stated before, there are several upper bounds of $\sigma(n)/n$ in the literature. All these estimates are of the equivalent form
\begin{equation}\label{eq6363.150h} 
	\frac{\sigma(n)}{n}\leq e^{\gamma} \log \log n+O\left(\frac{1}{(\log \log n)^2}\right),
\end{equation}
see \cite[Theorem 1.3]{AN2023},  et alii. Since
\begin{equation}\label{eq6363.150c} 
	e^{-c\sqrt{\log \log n}}\ll \frac{1}{(\log \log n)^b}
\end{equation}
for any constant $b\geq0$, the result in \hyperlink{thm1.1}{Theorem} \ref{thm1.1} seems to be the best possible short of the conditional result $\sigma(n)/n<e^{\gamma} \log \log n$.

\section{The Limits Infima} \label{s300}\hypertarget{s300}
The sum of divisors function is an oscillatory function, its values oscillate from its minimum $\sigma(n)=n+ 1$ at the prime integers $n$ to its maximum $\sigma(n)< 2n \log \log n$ at the extremely abundant integers $n$. The suprema of the sum of divisors function over various subsets of integers are known.
\begin{thm} \label{thm300.1}\hypertarget{thm300.1} {\normalfont (\cite{GT13})} Let $n \geq 1$ be an integer. Then
	\begin{enumerate} [font=\normalfont, label=(\roman*)]
		\item The limit supremum over the integers is
		$$
		\lim_{n \to \infty} \sup \frac{\sigma(n)}{n \log \log n}=e^{\gamma}.
		$$
		\item The limit supremum over the odd integers is
		$$
		\lim_{\text{odd }n \to \infty} \sup \frac{\sigma(n)}{n \log \log n}=\frac{e^{\gamma}}{2}.$$
		
		\item The limit supremum over the squarefree integers is
		$$
		\lim_{\text{squarefree }n \to \infty} \sup \frac{\sigma(n)}{n \log \log n}=\frac{6e^{\gamma}}{\pi^2}.$$
	\end{enumerate}
\end{thm}

Refer to \cite[p.\ 353]{HW79}, and similar references for proofs and other details on the maximal order of this function. \\


\section{Other Result} \label{S700}
Another technique for proving the sum of divisors inequality for the sequence of integers $n_r=n^r$, with $r  \geq 2$ is provided in this section. The verification involves the assumption that the inequality is false for all colossally abundant integers to derive a reductio ad absurdum. 

\begin{thm} \label{thm700.1}  Let $n \geq 5041$ be a sufficiently large integer, and let $r \geq 2$. Then
\begin{equation} \label{108}
 \sigma(n^r)< e^{\gamma} n \log \log n.
\end{equation}
\end{thm}

\begin{proof}[\textbf{Proof}] Let $n \geq 5041$ be a large colossally abundant integer, and let $r \geq 2$. Now, suppose that 
\begin{equation}\label{eq700.9}
 e^{\gamma}n^r \log \log n^r\leq\sigma(n^r)
\end{equation}
for all $r \geq 2$\\

\textbf{I.} Applying \hyperlink{lem400.1}{Lemma} \ref{lem400.1} and  \hyperlink{lem500.1}{Lemma} \ref{lem500.1} in succession reduce the right side to
\begin{eqnarray} \label{eq700.6}
\frac{\sigma(n^r)}{n^r} &<& \prod_{p \mid n^r} 
\left (1 -\frac{1}{p} \right )^{-1}  \\
&=& \prod_{p \mid n} 
\left (1 -\frac{1}{p} \right )^{-1} \nonumber \\
&\leq & \prod_{p \leq x} 
\left (1 -\frac{1}{p} \right )^{-1} \nonumber \\
&<& e^{\gamma} \log x \left (1 +\frac{1}{2(\log x)^2} \right ) \nonumber \\
&=& e^{\gamma} \log \log n \left (1 +\frac{1}{2(\log \log n)^2} \right )+O\left ( \frac{1}{(\log \log n)^2} \right ) \nonumber,
\end{eqnarray}
where $x=(\log n)(1+O(1/(\log \log n)^2))$, this follows from Lemma \ref{lem500.1}, and routine calculations. \\

\textbf{II.} By the hypothesis in equation (\ref{eq700.9}), and elementary calculations, the left side reduces to
\begin{eqnarray} \label{eq700.7}
\frac{\sigma(n^r)}{n^r} &\geq& e^{\gamma} \log \log n^r  \\
&=&e^{\gamma} \log \log n + e^{\gamma} \log \log r\nonumber .
\end{eqnarray}

\textbf{III.} Combining equations (\ref{eq700.6}) and  (\ref{eq700.7}) yield
\begin{eqnarray} \label{eq700.009a}
e^{\gamma} \log \log n \left (1 +\frac{1}{2(\log \log n)^2} \right )+O\left ( \frac{1}{(\log \log n)^2} \right )  &>& \frac{\sigma(n^r)}{n^r}  \\
&\geq&e^{\gamma} \log \log n + e^{\gamma} \log \log r   \nonumber.
\end{eqnarray}

Re expressing (\ref{eq700.9b}) in the equivalent form
\begin{equation}
1 +\frac{1}{2(\log \log n)^2} +O\left ( \frac{1}{(\log \log n)^3} \right )>1+\frac{\log r}{\log \log n}  
\end{equation}
shows that the left side converges to $1$ at a faster rate than the right side by a factor of $1/\log \log n$. Clearly, there is a contradiction for all $r \geq 2$. Ergo, $\sigma(n^r)< e^{\gamma}n \log \log n$ for all sufficiently large integers $n \geq 1$, and a fixed parameter $r \geq 2$.
\end{proof}



\begin{thebibliography}{999}
\bibitem{AN2023}Axler, C.; Nicolas, J. \textit{\color{red}Large values of $n/\varphi(n)$ and $\sigma(n)/n$.} Acta Arithmetica, First Online Preprint, September, 2023.

\bibitem{AC17} Axler, Christian. \textit{\color{red}A new upper bound for the sum of divisors function.} Bull. Aust. Math. Soc. 96 (2017), no. 3, 374-379.

\bibitem{AE44} L. Alaoglu and P. Erdos, \textit{\color{red}On highly composite and similar numbers.} Trans. Amer: Math. Soc. 56 (1944), 448-469.

\bibitem{AJ07} Akbary A., Friggstad Z. and Juricevic R., \textit{\color{red}Explicit upper bounds for $\prod_{x \leq x}p/(p-1)$.} Contrib. Discrete Math. 2(2) (2007), 153--160. 

\bibitem{BN07} William D. Banks, Derrick N. Hart, Pieter Moree, C. Wesley Nevans, \textit{\color{red}The Nicolas and Robin inequalities
with sums of two squares.} Monatsh. Math. 157, no. 4,(2009) 303-322. \url{http://arxiv.org/abs/0710.2424}.

\bibitem{BK06} Briggs, K. \textit{\color{red}Abundant numbers and the Riemann hypothesis.} Experiment. Math. 15 (2006), no. 2, 251-256.

\bibitem{BR2018} Michael A. Bennett, Greg Martin, Kevin O'Bryant, Andrew Rechnitzer. \textit{\color{red}Explicit bounds for primes in arithmetic progressions}. http://arxiv.org/abs/1802.00085.

\bibitem{BT15} Broughan K. and Trudgian T., \textit{\color{red}Robins inequality for 11-free integers.} Integers 15 (2015), Article ID A12, 5 pages.


\bibitem{CP1985} N. Costa Pereira. \textit{\color{red} Estimates for the Chebyshev function $\theta(x)$ and $\psi(x)$}. Math. Comp., 44(169):211--221, 1985.


\bibitem{CS06} Y.-J. Choie, N. Lichiardopol, P. Moree, P. Sole, \textit{\color{red}On Robin's criterion for the Riemann Hypothesis.} J. Theor. Nombres Bordeaux 19 (2007), no. 2, 357-372.
\url{http://arxiv.org/abs/math/0604314}.


\bibitem{DH1980} Davenport, Harold. \textit{\color{red}Multiplicative number theory.} Springer-Verlag Berlin, Heidelberg, New York, 1980. 

\bibitem{DP2010} P. Dusart. \textit{\color{red}Estimates of some functions over primes without R.H.} http://arxiv.org/abs/1002.0442.

\bibitem{DP2009} H. G. Diamond and J. Pintz. \textit{\color{red}Oscillation of Mertens product formula}. J. Theor. Nombres Bordeaux 21 (2009), no. 3, 523-533.


\bibitem{EE1985} Ellison, William; Ellison, Fern, Prime numbers. Wiley-Interscience Publication. New York; Hermann, Paris, 1985.

\bibitem{EL1888} Euler, Leonhard. Introduction to analysis of the infinite. Book I. Translated from the Latin and with an introduction by John D. Blanton. 


\bibitem{FK2019} Ford, K.  \textit{\color{red}Zero-free regions for the Riemann zeta function.} 
\url{http://arxiv.org/abs/1910.08205}. 

\bibitem{FS03}Finch,  S.R. \textit{\color{red}Mathematical constants.} Encyclopedia of Mathematics and its Applications 94, (Cambridge
University Press, Cambridge, 2003).

\bibitem{GR07} A. Grytczuk, \textit{\color{red}Upper bound for sum of divisors function and Riemann Hypothesis.} Tsukuba J. Math, Vol. 31 No.1,(2007) 67-75.

\bibitem{GT13} Gronwall, T. H. \textit{\color{red}Some asymptotic expressions in the theory of numbers.} Trans. Amer. Math. Soc. 14 (1913), no. 1, 113-122.

\bibitem{IA2003} Ivic, Aleksandar. \textit{\color{red}The Riemann zeta-function. Theory and applications}. Wiley, New York; Dover Publications, Inc., Mineola, NY, 2003.



\bibitem{IA1942} A.E. Ingham, On two conjectures in the theory of numbers. Amer. J. Math., 64 (1942), 313-319.

\bibitem{HA16} Alexander Hertlein, \textit{\color{red}Robin's inequality for new families of integers.} \url{http://arxiv.org/abs/1612.05186}.

\bibitem{HW79} G. H. Hardy and E. M. Wright, \textit{\color{red}An Introduction to the Theory of Numbers.} 5th ed., Oxford University
Press, Oxford, 1979.

\bibitem{LJ00} Lagarias, J. C. \textit{\color{red}An Elementary Problem Equivalent to the Riemann Hypothesis.} \url{http://arxiv.org/abs/math/0008177}.

\bibitem{LJ13} Lagarias, J. C. \textit{\color{red}Euler's constant: Euler's work and modern developments.} \url{http://arxiv.org/abs/1303.1856}.

\bibitem{LJ2021} Lichtman, J. \textit{\color{red}Averages of the Mobius function on shifted primes}. http://arxiv.org/abs/2009.08969.

\bibitem{LJ1963} Linnik, J. V. \textit{\color{red}The dispersion method in binary additive problems.} (Leningrad 1961), Transl. Math. monographs, Vol. 4, Amer. Math. Soc., Providence, R. I., 1963; Chapter 8.

\bibitem{LJ2015} Lay, J. P. S. \textit{\color{red}Sign changes in Mertens first and second theorems}, http://arxiv.org/abs/1505.03589.

\bibitem{LY2014} Lamzouri, Y. \textit{\color{red}A bias in Mertens product formula}, http://arxiv.org/abs/1410.3777v2. 


\bibitem{LZ2007B} Languasco, A.; Zaccagnini, A. \textit{\color{red}On the constant in the Mertens product for arithmetic progressions. I}. http://arxiv.org/abs/0706.2807. 

\bibitem{LZ2007} Languasco, A.; Zaccagnini, A. \textit{\color{red}A note on Mertens' formula for arithmetic progressions.} J. Number Theory  127  (2007),  no. 1, 37-46. 

\bibitem{LN2009} Lidl, R.; Niederreiter, H. \textit{\color{red}Finite fields.} With a foreword by P. M. Cohn. Second edition. Encyclopedia of Mathematics and its Applications, 20. Cambridge University Press, Cambridge, 1997.

\bibitem{MV2007} Montgomery, H. L.; Vaughan, R. C. \textit{\color{red}Multiplicative number theory. I. Classical theory}. Cambridge University Press, Cambridge, 2007.

\bibitem{KP2020} Kinlaw, P.; Kobayashi, M.; Pomerance, C. \textit{\color{red}On the equation $\varphi(n)=\varphi(n+1)$}. Acta Arithmetica 196 (2020), 69-92.

\bibitem{NJ83} J. L. Nicolas, \textit{\color{red}Petites valeurs de la fonction d Euler.} J. Number Theory 17 (1983) 375-388.

\bibitem{NW2000} Narkiewicz, W. \textit{\color{red}The development of prime number theory. From Euclid to Hardy and Littlewood}. Springer Monographs in Mathematics. Springer-Verlag, Berlin, 2000.

\bibitem{NS12} S. Nazardonyavi, \textit{\color{red}Superabundant numbers, their subsequences and the Riemann.} \url{http://arxiv.org/abs/1211.2147}.	 
	 
\bibitem{OI2026} OpenAI. Proof of the Quasi-Riemann Hypothesis. Preprint, October 6, 2026.

\bibitem{QRH2026} Mathoverflow. Quasi Riemann Hypothesis, 7/8Multiplicative functions in short interval. \href{https://mathoverflow.net/questions/515968/consequences-of-quasi-rh}{Quasi RH.}




\bibitem{RG84} Guy Robin, \textit{\color{red}Grandes valeurs de la fonction somme des diviseurs et hypothese de Riemann}. J. Math. Pures Appl. (9) 63 (1984), 187-213.

\bibitem{RP1996} Ribenboim, Paulo. \textit{\color{red}The new book of prime number records.} Berlin, New York: Springer-Verlag, 1996.

\bibitem{RS1962} Rosser, J.B.; Schoenfeld, L. \textit{\color{red}Approximate formulas for some functions of prime numbers}. Illinois J. Math. 6 (1962) 64--94.


\bibitem{RS1915}Ramanujan, S. \textit{\color{red}Highly composite numbers}. Proc. London Math. Soc. (2) 14 (1915), 347-409.

\bibitem{RS1997}Ramanujan, S. \textit{\color{red}Highly composite numbers.} Ramanujan J. 1 (1997), 119-153.

\bibitem{SL1976} Schoenfeld, Lowell. \textit{\color{red}Sharper bounds for the Chebyshev functions $\theta (x)$ and $\psi (x)$. II.} Math. Comp. 30 (1976), no. 134, 337-360.



\bibitem{SP10} Patrick Sole, Michel Planat, \textit{\color{red}Robin inequality for 7-free integers.} \url{http://arxiv.org/abs/1012.0671}.


\bibitem{VM2005} Mark B. Villarino, \textit{\color{red}Mertens' Proof of Mertens' Theorem.} \url{http://arxiv.org/abs/math/0504289}.



\bibitem{WA2013} Weingartner, A. \textit{\color{red}The distribution functions of  $\sigma(n)/n$  and  $n/\varphi(n)$, II.} J. Number Theory 132 (2012), no. 12, 2907--2921. \href{https://mathscinet.ams.org/mathscinet-getitem?mr=MR2965199}{MR2965199}


\bibitem{WM07} Marek Wojtowicz, \textit{\color{red}Robins inequality and the Riemann hypothesis.} Proc. Japan Acad. Ser. A Math. Sci. Volume 83, Number 4 (2007), 47-49.



 
\end{thebibliography}
\end{document}